# On vanishing diffusivity selection for the advection equation

Giulia Mescolini[1] · Jules Pitcho[2,3] · Massimo Sorella[4]

**Abstract**
We study the advection equation along vector fields singular at the initial time. More precisely, we prove that for divergence-free vector fields in $L^1_{loc}((0, T]; BV(\mathbb{T}^d; \mathbb{R}^d)) \cap L^2((0, T) \times \mathbb{T}^d; \mathbb{R}^d))$, there exists a unique vanishing diffusivity solution. This class includes the vector field constructed by Depauw in [13], for which there are infinitely many distinct bounded solutions to the advection equation.

**Mathematics Subject Classification** 35A02 · 35D30 · 35Q49 · 34A12

## 1 Introduction

We here study the initial value problem for the advection equation

$$\begin{cases} \partial_t \rho + \boldsymbol{b} \cdot \nabla \rho = 0, \\ \rho(0, x) = \rho_{in}(x), \end{cases} \quad \text{(PDE)}$$

posed on $[0, T] \times \mathbb{T}^d$ where $\boldsymbol{b} = \boldsymbol{b}(t, x)$ is a given vector field, $\rho = \rho(t, x)$ is an unknown real-valued function, $\rho_{in} = \rho_{in}(x)$ is the bounded initial datum and $\mathbb{T}^d \cong \mathbb{R}^d/\mathbb{Z}^d$ is the periodic torus. In this work, we consider vector fields, which are divergence-free in the sense of distributions, and we study weak solution of (PDE).

**Definition 1.1** Consider an integrable vector field $\boldsymbol{b} : [0, T] \times \mathbb{T}^d \to \mathbb{R}^d$, and an initial datum $\rho_{in} \in L^\infty(\mathbb{T}^d)$. We shall say that $\rho \in L^\infty((0, T) \times \mathbb{T}^d)$ is a bounded weak solution of (PDE) along $\boldsymbol{b}$, if for every $\phi \in C_c^\infty([0, T) \times \mathbb{T}^d)$, we have

$$\int_0^T \int_{\mathbb{T}^d} \left[ \rho(t, x) \partial_t \phi(t, x) + \rho(t, x) \boldsymbol{b}(t, x) \cdot \nabla \phi(t, x) \right] dx dt = -\int_{\mathbb{T}^d} \rho_{in}(x) \phi(0, x) dx. \quad (1.1)$$

✉ Jules Pitcho
  jules.pitcho@gssi.it

  Giulia Mescolini
  giulia.mescolini@epfl.ch

  Massimo Sorella
  m.sorella@imperial.ac.uk

[1] Institute of Mathematics, EPFL, Station 8, 1015 Lausanne, Switzerland

[2] ENS de Lyon, UMPA, 46 allée d'Italie, 69364 Lyon, France

[3] GSSI, Via Michele Iacobucci, 2, 67100 L'Aquila, Italy

[4] Department of Mathematics, Imperial College of London, London SW7 2AZ, UK

A regularisation and compactness argument shows the existence of bounded weak solutions if we assume $\int_0^T \|[\text{div } \boldsymbol{b}(s,\cdot)]^-\|_{L^\infty(\mathbb{T}^d)} ds < +\infty$. By the Cauchy-Lipschitz theory, uniqueness also holds when the time integral of the spatial Lipschitz constant of the vector field is finite. However uniqueness of bounded weak solutions can fail when the vector field is rough. One example is given by Depauw in [13], where he constructs a bounded, divergence-free vector field, for which there are infinitely many distinct bounded weak solutions to (PDE). We introduce the diffusivity parameter $\nu \in (0,1)$ and we study weak solutions to the advection–diffusion equation

$$\begin{cases} \partial_t \rho^\nu + \text{div}(\boldsymbol{b}\rho^\nu) - \nu\Delta\rho^\nu = 0, \\ \rho(0,\cdot) = \rho_{in}. \end{cases} \quad (\nu\text{-PDE})$$

for which uniqueness of bounded weak solutions is restored assuming only that $\boldsymbol{b} \in L^2((0,T) \times \mathbb{T}^d; \mathbb{R}^d)$ see [5, Theorem 2.7 and Theorem 3.3] and Sect. 2. In this context, we introduce vanishing diffusivity solutions as a subclass of weak solutions.

**Definition 1.2** Consider a vector field $\boldsymbol{b} \in L^2((0,T) \times \mathbb{T}^d)$, and an initial datum $\rho_{in} \in L^\infty(\mathbb{T}^d)$. For every $\nu \in (0,1)$, consider the unique bounded solution $\rho^\nu$ of ($\nu$-PDE) along $\boldsymbol{b}$ with initial datum $\rho_{in}$ independent of $\nu$. We shall say that $\rho \in L^\infty((0,T) \times \mathbb{T}^d)$ is a vanishing diffusivity solution along $\boldsymbol{b}$ with initial datum $\rho_{in}$, if there exists a sequence of real numbers $(\nu_i)_{i\in\mathbb{N}}$ in $(0,1)$ such that $\nu_i \downarrow 0$ as $i \to +\infty$, and $\rho^{\nu_i}$ converges to $\rho$ weakly-star in $L^\infty((0,T) \times \mathbb{T}^d)$ as $i \to +\infty$.

We here study uniqueness of vanishing diffusivity solutions.

**Question 1.3** *What are the minimal assumptions on a divergence-free vector field $\boldsymbol{b} \in L^2((0,T) \times \mathbb{T}^d; \mathbb{R}^d)$ so that there exists a unique vanishing diffusivity solution for every bounded initial datum?*

This question is natural. For scalar conservation laws with smooth flux in one-space dimension, Kružkov [18] proves uniqueness and continuous dependance with respect to initial data of vanishing viscosity solutions, even though uniqueness of weak solutions fails. Bressan, Guerra and Shen [7] also prove uniqueness and continuous dependance with respect to initial data of vanishing viscosity solutions for scalar conservation laws in one-space dimension with possibly discontinuous flux satisfying an additional condition. For strictly hyperbolic systems in one-space dimension, Bianchini and Bressan [3] prove uniqueness and continuous dependance with respect to initial data of vanishing viscosity solutions. Note finally that in several space dimensions, hyperbolic systems can be ill-posed as the example [6] shows.

The problem of selection via vanishing diffusivity has also been studied for the transport equation. In [1] Ambrosio, building upon the work of DiPerna and Lions [14], proves that for divergence-free vector fields $\boldsymbol{b} \in L^1((0,T); BV(\mathbb{T}^d; \mathbb{R}^d))$ there exists a unique bounded weak solution to the transport equation. This implies directly uniqueness of the vanishing diffusivity solution. Vanishing diffusivity is also shown to be a selection criterion for more general initial data by Bonicatto, Ciampa and Crippa in [4]. In [10] Colombo, Crippa and the third author construct a $C^\alpha((0,2) \times \mathbb{T}^d)$ divergence-free vector field, with $\alpha \in (0,1)$ fixed but arbitrary, for which there are at least two distinct vanishing diffusivity solutions. In this example the vector field has a single time of singularity at $t = 1$. In [16] Huysmans and Titi construct a bounded, divergence–free vector field for which there exists a unique *inadmissible* vanishing diffusivity solution, in the sense that the $L^2$ norm of this solution is not monotonous. A related question is anomalous dissipation. We say that a vector field $\boldsymbol{b} \in L^2((0, T) \times \mathbb{T}^d; \mathbb{R}^d)$ has anomalous dissipation for an initial datum $\rho_{in}$ if the family of unique weak solutions $\{\rho^\nu\}_{\nu \in (0,1)}$ of ($\nu$ -PDE) satisfies

$$\limsup_{\nu \to 0} \nu \int_0^T \int_{\mathbb{T}^d} |\nabla \rho^\nu|^2 > 0.$$

In [22], it is shown that anomalous dissipation implies non–uniqueness of the backward advection equation. But for vector fields which have anomalous dissipation, there may be a unique vanishing diffusivity solution as the construction of [15] shows. See also [2, 8, 17] for works on anomalous dissipation.

An alternative to studying the vanishing diffusivity limit is to regularise the vector field, and to study the limit of vanishing regularisation. A series of works [9–11, 19] study the absence of selection by such schemes. The second author of this paper further remarks in [21] that for bounded, divergence-free vector fields in $L^1_{loc}((0, T]; BV(\mathbb{T}^d; \mathbb{R}^d))$ regularisation by mollification selects a unique solution of (PDE) for any integrable initial datum. He also presents a Lagrangian counterpart to this selection result in [20], where interestingly he shows that for Depauw's example [13], a probability measure concentrated on several distinct integral curves is selected for Lebesgue almost every initial datum. In this paper, we apply the ideas developed by the second author in [21] to show that the vanishing diffusivity is a selection criterion for the advection equation along a divergence-free Borel vector field $\boldsymbol{b} : [0, T] \times \mathbb{T}^d \to \mathbb{R}^d$ in $L^1_{loc}((0, T]; BV(\mathbb{T}^d; \mathbb{R}^d)) \cap L^2((0, T) \times \mathbb{T}^d; \mathbb{R}^d)$.

Recall that a *standard mollifier* is a function $w \in C_c^\infty((0, 1)^d)$ such that $w \geq 0$, and $\int_{(0,1)^d} w(x)dx = 1$, and that it induces a mollification family $\{w^\delta\}_{\delta \in (0,1)}$ on $\mathbb{T}^d$ where we denote $w^\delta(x) = w(x/\delta)/\delta^d$, and we identify $w^\delta$ with its periodisation to a function on $\mathbb{T}^d$. Recall also that the mollification of $\boldsymbol{b}$ with $w^\delta$ is defined for every $x \in \mathbb{T}^d$ and almost every $t \in [0, T]$ as

$$(\boldsymbol{b} * w^\delta)(t, x) = \int_{\mathbb{T}^d} \boldsymbol{b}(t, y) w^\delta(x - y) dy.$$

Moreover by Young's convolution inequality, we have

$$\|\nabla(\boldsymbol{b} * w^\delta)\|_{L^1((0,T);L^\infty(\mathbb{T}^d))} \leq \|\boldsymbol{b}\|_{L^1((0,T);L^1(\mathbb{T}^d))} \|\nabla w^\delta\|_{L^\infty(\mathbb{T}^d)}.$$

Therefore by the Cauchy-Lipschitz theory, for any initial datum $\rho_{in} \in L^\infty(\mathbb{T}^d)$, there exists a unique bounded weak solution of (PDE) along $\boldsymbol{b} * w^\delta$. We now state our result.

**Theorem 1.4** *Consider a divergence-free vector field $\boldsymbol{b} \in L^1_{loc}((0, T]; BV(\mathbb{T}^d; \mathbb{R}^d)) \cap L^2((0, T) \times \mathbb{T}^d; \mathbb{R}^d))$, and an initial datum $\rho_{in} \in L^\infty(\mathbb{T}^d)$. Then there exists a unique vanishing diffusivity solution of (PDE) along $\boldsymbol{b}$ with initial datum $\rho_{in}$. Futhermore,*

*(i) for every standard mollifier $w$, the unique bounded weak solutions $\tilde{\rho}^\delta$ of (PDE) along $\boldsymbol{b} * w^\delta$ converge weakly-star in $L^\infty((0, T) \times \mathbb{T}^d)$ to the vanishing diffusivity solution as $\delta \downarrow 0$;*

*(ii) there is no anomalous dissipation, meaning that the unique bounded weak solutions $\rho^\nu$ of ($\nu$ -PDE) along $\boldsymbol{b}$ with initial datum $\rho_{in}$ satisfy*

$$\limsup_{\nu \downarrow 0} \nu \int_0^T \int_{\mathbb{T}^d} |\nabla \rho^\nu(t, x)|^2 dx dt = 0.$$

**Remark 1.5** In contrast with the selection result in [21], we do not have the optimal exponents between the vector field and the initial datum which insure that weak solutions are well-defined. This is because, for divergence-free vector fields in $L^2((0,T) \times \mathbb{T}^d; \mathbb{R}^d)$, we can only prove uniqueness for bounded weak solutions of ($\nu$ -PDE).

**Remark 1.6** It would be interesting to determine whether the vanishing diffusivity solutions of Theorem 1.4 depend continuously on the initial data.

### 1.1 Our argument and organisation of the paper

For the selection part of our result, we relate forwards-in-time vanishing diffusivity selection to a backwards-in-time uniqueness of vanishing diffusivity solution by duality. The backwards-in-time uniqueness for vanishing diffusivity solutions follows because of uniqueness of bounded solution for strictly positive times as well as as uniqueness of the trace up to time $t = 0$ in an $L^\infty$ weak-star sense. This strategy draws from the work [21] of the second author. The proof of the absence of anomalous dissipation property draws from [12].

In Sect. 2, we gather some useful existence and uniqueness results for weak solutions of (PDE) and ($\nu$ -PDE). In Sect. 3, we thoroughly analyse the backwards problems and prove a uniqueness result for vanishing diffusivity solutions, as well as a duality formula. In Sect. 4, we conclude the proof of our main theorem.

## 2 Preliminaries

The aim of this section is to state some known results we will need for the proof of our main theorem.

### 2.1 The non-diffusive problem

We recall the following well-posedness result on (PDE) from [1] in the context of divergence-free vector fields.

**Theorem 2.1** *Consider a divergence-free vector field $\boldsymbol{b} \in L^1((0,T); BV(\mathbb{T}^d; \mathbb{R}^d))$, and an initial datum $\rho_{in} \in L^\infty(\mathbb{T}^d)$. Then there exists a unique bounded weak solution of (PDE) along $\boldsymbol{b}$ with initial datum $\rho_{in}$. Furthermore, up to redefining $\rho$ on a zero measure set we have $\rho \in C([0,T]; L^2(\mathbb{T}^d))$, and*

$$\|\rho(t,\cdot)\|_{L^2(\mathbb{T}^d)} = \|\rho_{in}\|_{L^2(\mathbb{T}^d)} \quad \text{for every } t \in [0,T].$$

### 2.2 The diffusive problem

The following remark summarizes the solution theory for ($\nu$ -PDE) along a smooth vector field and with a smooth initial datum.

**Remark 2.2** Let $(\Omega, (\mathcal{F}_t)_{t \geq 0}, \mathbb{P})$ be a filtered probability space. For $t \in [0, T]$, consider the following backward stochastic differential equation

$$dX^\nu_{t,s} = \boldsymbol{b}(s, X^\nu_{t,s})ds + \sqrt{2\nu}d\boldsymbol{W}_s \quad \text{with} \quad X^\nu_{t,t}(x) = x, \tag{2.1}$$

where $W$ is a $\mathbb{T}^d$-valued Brownian motion adapted to the backward filtration such that we have $W_t = 0$ almost surely. Assume now that the vector field $\boldsymbol{b}$ is divergence-free and smooth. Then for almost every $\omega \in \Omega$, we have $X_{t,s}(\omega)$ is incompressible, namely $X_{t,s}(\omega, \cdot)_{\#}\mathscr{L}^d = \mathscr{L}^d$. We also recall that the Feynman-Kac formula, which gives the following expression for solutions of ($\nu$ -PDE) with smooth initial datum

$$\rho^\nu(t, x) = \mathbb{E}[\rho_{in}(X^\nu_{t,0}(x))].$$

It can be checked that $\rho^\nu \in C([0, T]; L^2(\mathbb{T}^d))$, and a computation shows that

$$\sup_{t \in [0,T]} \|\rho^\nu(t, \cdot)\|_{L^\infty(\mathbb{T}^d)} \leq \|\rho_{in}\|_{L^\infty(\mathbb{T}^d)}.$$

Let us now define bounded weak solutions of the diffuse problem.

**Definition 2.3** Consider a vector field $\boldsymbol{b} \in L^1((0, T) \times \mathbb{T}^d; \mathbb{R}^d)$ and an initial datum $\rho_{in} \in L^\infty(\mathbb{T}^d)$. We shall say that $\rho \in L^\infty((0, T) \times \mathbb{T}^d)$ is a bounded weak solution of ($\nu$ -PDE) along $\boldsymbol{b}$ with initial datum $\rho_{in}$, if for every $\phi \in C^\infty_c([0, T) \times \mathbb{T}^d)$, we have

$$\int_0^T \int_{\mathbb{T}^d} \Big[\rho(t, x)\partial_t\phi(t, x) + \boldsymbol{b}(t, x)\rho(t, x) \cdot \nabla\phi(t, x) - \nu\rho(t, x)\Delta\phi(t, x)\Big]dxdt$$
$$= -\int_{\mathbb{T}^d} \rho_{in}(x)\phi(0, x)dx. \tag{2.2}$$

The following well-posedness theorem will be useful. We also provide a proof.

**Theorem 2.4** *Consider a divergence-free vector field $\boldsymbol{b} \in L^2((0, T) \times \mathbb{T}^d; \mathbb{R}^d)$, an initial datum $\rho_{in} \in L^\infty(\mathbb{T}^d)$, and a real number $\nu > 0$. Then there exists a unique weak solution $\rho^\nu \in L^\infty((0, T) \times \mathbb{T}^d) \cap L^2((0, T); H^1(\mathbb{T}^d))$ of ($\nu$ -PDE) along $\boldsymbol{b}$ with initial datum $\rho_{in}$ such that*

$$\|\rho^\nu\|_{L^\infty((0,T)\times\mathbb{T}^d)} \leq \|\rho_{in}\|_{L^\infty(\mathbb{T}^d)}.$$

*Moreover $\rho^\nu$ belongs to $C([0, T]; L^2(\mathbb{T}^d))$ up to modification on a zero measure set, and for every $t \in [0, T]$*

$$\|\rho^\nu(t, \cdot)\|^2_{L^2(\mathbb{T}^d)} + 2\nu\|\nabla\rho^\nu\|^2_{L^2((0,t)\times\mathbb{T}^d)} \leq \|\rho_{in}\|^2_{L^2(\mathbb{T}^d)}. \tag{2.3}$$

***Proof*** **Step 1** (Existence) Let $\{w^\delta\}_{\delta>0}$ be a standard family of mollifiers; for example, consider $w^\delta(x) = \delta^{-d}w(x/\delta)$ with $w$ smooth and compactly supported. Define $\boldsymbol{b}^\delta = \boldsymbol{b} * w^\delta$ and $\rho^\delta_{in} = \rho_{in} * w^\delta$. We consider the problem with mollified data

$$\begin{cases} \partial_t\rho^{\nu,\delta} + \text{div}(\boldsymbol{b}^\delta\rho^{\nu,\delta}) - \nu\Delta\rho^{\nu,\delta} = 0, \\ \rho^{\nu,\delta}(0, \cdot) = \rho^\delta_{in}. \end{cases} \tag{2.4}$$

By the classical theory detailed in Remark 2.2, (2.4) admits a smooth solution that we can bound in $L^2((0, T); H^1(\mathbb{T}^d))$: indeed, multiplying by $\rho^{\nu,\delta}$ and integrating by parts both in space and time, we get

$$\int_0^t \partial_s \int_{\mathbb{T}^d} \rho^{\nu,\delta}(s, x)\rho^{\nu,\delta}(s, x)dxds - 2\nu \int_0^t \int_{\mathbb{T}^d} \Delta\rho^{\nu,\delta}(s, x)\rho^{\nu,\delta}(s, x)dxds = 0,$$

which implies for a.e. $t \in [0, T]$

$$\|\rho^{\nu,\delta}(t, \cdot)\|^2_{L^2(\mathbb{T}^d)} + 2\nu\|\nabla\rho^{\nu,\delta}\|^2_{L^2((0,t)\times\mathbb{T}^d)} \leq \|\rho^\delta_{in}\|^2_{L^2(\mathbb{T}^d)} \leq \|\rho_{in}\|^2_{L^2(\mathbb{T}^d)}. \tag{2.5}$$

Together with Remark 2.2, this gives the uniform bounds in $\delta$:

$$\begin{aligned}\|\rho^{\nu,\delta}\|_{L^\infty((0,T)\times\mathbb{T}^d)} &\leq \|\rho_{in}\|_{L^\infty(\mathbb{T}^d)},\\ \|\rho^{\nu,\delta}\|_{L^2((0,T);H^1(\mathbb{T}^d))} &\leq C_{T,\nu}\|\rho_{in}\|_{L^2(\mathbb{T}^d)},\end{aligned} \quad (2.6)$$

where the constant $C_{T,\nu}$ depends on $T$ and $\nu$ only. Hence, the family $\{\rho^{\nu,\delta}\}_{\delta>0}$ is uniformly bounded and we can extract a subsequence which converges weakly-star in $L^\infty((0,T)\times\mathbb{T}^d)$, and in $L^2((0,T);H^1(\mathbb{T}^d))$ to a function $\rho^\nu$, which is a weak solution of ($\nu$-PDE) since the equation is linear and by weak lower semicontinuity, $\rho^\nu$ satisfies the bounds

$$\begin{aligned}\|\rho^\nu\|_{L^\infty((0,T)\times\mathbb{T}^d)} &\leq \|\rho_{in}\|_{L^\infty(\mathbb{T}^d)},\\ \|\rho^\nu\|_{L^2((0,T);H^1(\mathbb{T}^d))} &\leq C_{T,\nu}\|\rho_{in}\|_{L^2(\mathbb{T}^d)}.\end{aligned} \quad (2.7)$$

**Step 2** (Uniqueness) Let $v$ be a solution of ($\nu$-PDE) in $L^\infty((0,T)\times\mathbb{T}^d) \cap C([0,T];L^2(\mathbb{T}^d))$ with $\rho_{in} = 0$ (which corresponds to the problem solved by the difference of two solutions), and note that the mollification of the solution $v^\delta := v * w^\delta$ satisfies:

$$\begin{cases}\partial_t v^\delta - \nu \Delta v^\delta + \text{div}(\boldsymbol{b}v^\delta) = r^\delta,\\ v^\delta(0,\cdot) = 0,\end{cases} \quad (2.8)$$

with $r^\delta := \text{div}(\boldsymbol{b}v^\delta - (\boldsymbol{b}v) * w^\delta)$. We can prove that $r^\delta$ is bounded in $L^2((0,T);\dot{H}^{-1}(\mathbb{T}^d))$ and

$$\|r^\delta\|_{L^2((0,T);\dot{H}^{-1}(\mathbb{T}^d))} \to 0 \text{ as } \delta \downarrow 0.$$

Indeed, for any $\phi \in L^2((0,T);\dot{H}^1(\mathbb{T}^d))$, we have:

$$\left|\int_0^T \int_{\mathbb{T}^d} r^\delta(s,x)\phi(s,x)dxds\right|$$
$$= \left|\int_0^T \int_{\mathbb{T}^d} (\boldsymbol{b}(s,x)v^\delta(s,x) - (\boldsymbol{b}v) * w^\delta(s,x)) \cdot \nabla\phi(s,x)dxds\right|$$
$$= \left|\int_0^T \int_{\mathbb{T}^d} \left(\int_{\mathbb{T}^d} (\boldsymbol{b}(s,x) - \boldsymbol{b}(s,x-y))v(s,x-y)w^\delta(y)dy\right) \cdot \nabla\phi(s,x)dxds\right|$$
$$\leq \int_K |w(z)| \int_0^T \int_{\mathbb{T}^d} |\boldsymbol{b}(s,x) - \boldsymbol{b}(s,x-\delta z)||v(s,x-\delta z)||\nabla\phi(s,x)|dxdsdz,$$

using the change of variables $z = y/\delta$ and recalling that $\text{supp}(w) \subset K$ with $K \subset \mathbb{T}^d$ compact. By dominated convergence theorem, observing that the integrand converges in measure to 0, we can prove that:

$$\lim_{\delta\to 0}\int_0^T \int_{\mathbb{T}^d} |\boldsymbol{b}(s,x) - \boldsymbol{b}(s,x-\delta z)||v(s,x-\delta z)||\nabla\phi(s,x)|dxds = 0. \quad (2.9)$$

Now, let $t \in [0,T]$, and multiply (2.8) by $v^\delta$ and integrating by parts in time and space over $(0,t)\times\mathbb{T}^d$, we get:

$$\frac{\|v^\delta(t,\cdot)\|^2_{L^2(\mathbb{T}^d)}}{2} + \nu\int_0^t \int_{\mathbb{T}^d} |\nabla v^\delta(s,x)|^2 dxds = \int_0^t \int_{\mathbb{T}^d} r^\delta(s,x)v^\delta(s,x)dxds, \quad (2.10)$$

The advection term vanishes, because

$$\int_0^t \int_{\mathbb{T}^d} \operatorname{div}(\boldsymbol{b}(s,x)v^\delta(s,x))v^\delta(s,x)dxds = -\int_0^t \int_{\mathbb{T}^d} (\boldsymbol{b}(s,x) \cdot \nabla v^\delta(s,x))v^\delta(s,x)dxds$$

$$= -\int_0^t \int_{\mathbb{T}^d} \operatorname{div}(\boldsymbol{b}(s,x)v^\delta(s,x))v^\delta(s,x)dxds,$$

because div $\boldsymbol{b} = 0$. By (2.10) and Young's inequality, we get

$$\frac{\|v^\delta(t,\cdot)\|^2_{L^2(\mathbb{T}^d)}}{2} + \nu\|v^\delta\|^2_{L^2((0,t);\dot{H}^1(\mathbb{T}^d))} \leq \|r^\delta\|_{L^2((0,t);\dot{H}^{-1}(\mathbb{T}^d))}\|v^\delta\|_{L^2((0,t);\dot{H}^1(\mathbb{T}^d))} \quad (2.11)$$

$$\leq \nu\|v^\delta\|^2_{L^2((0,t);\dot{H}^1(\mathbb{T}^d))} + \frac{1}{4\nu}\|r^\delta\|^2_{L^2((0,t);\dot{H}^{-1}(\mathbb{T}^d))}. \quad (2.12)$$

The desired uniqueness result comes by observing that

$$\frac{\|v^\delta(t,\cdot)\|^2_{L^2(\mathbb{T}^d)}}{2} \leq \frac{1}{4\nu}\|r^\delta\|^2_{L^2((0,t);\dot{H}^{-1}(\mathbb{T}^d))} \longrightarrow 0 \quad \text{as } \delta \downarrow 0,$$

and that $t$ was arbitrary in $[0, T]$.

**Step 3** (Time continuity) Let us begin by showing that, up to redefining $\rho$ on a set of zero measure, we have $\rho^\nu \in C([0, T]; w - L^2(\mathbb{T}^d))$. Let $\mathcal{N}$ be a dense, countable subset of $L^2(\mathbb{T}^d)$ consisting of smooth functions, which exists by separability, and let $\phi \in \mathcal{N}$. Let also $\psi \in C_c^\infty((0, T])$. Define

$$\Theta_\phi(t) := \int_{\mathbb{T}^d} \rho^\nu(t,x)\phi(x)dx,$$

and note that by the weak formulation of ($\nu$ -PDE), we have

$$\int_0^T \Theta_\phi \partial_\tau \psi(\tau) d\tau = -\int_{\mathbb{T}^d} \rho_{in}(x)\phi(x)dx - \int_0^T \int \boldsymbol{b}(\tau,x)\rho^\nu(\tau,x) \cdot \nabla\phi(x)\psi(\tau)dxd\tau$$

$$+ \int_0^T \int_{\mathbb{T}^d} \rho^\nu(\tau,x)\Delta\phi(x)\psi(\tau)dxd\tau.$$

Since $\boldsymbol{b} \in L^2((0,T) \times \mathbb{T}^d; \mathbb{R}^d)$, $\rho^\nu \in L^\infty((0,T) \times \mathbb{T}^d)$, and $\phi$ is smooth, we conclude that $\partial_t \Theta_\phi \in L^1((0,T))$, hence $\Theta_\phi$ is continuous up to a set $\mathcal{T}_\phi$ of measure zero. Now let

$$\mathcal{T} := \bigcup_{\phi \in \mathcal{N}} \mathcal{T}_\phi,$$

which also has zero measure. For every $\phi \in \mathcal{N}$, we redefine $\Theta_\phi$ on $\mathcal{T}$ to be continuous, and we observe that since $|\Theta_\phi(t)| \leq \|\rho^\nu\|_{L^\infty((0,T);L^2(\mathbb{T}^d))}\|\phi\|_{L^2(\mathbb{T}^d)}$ for all $t \in [0, T]$, the bounded linear functional $L_t : \mathcal{N} \ni \phi \mapsto \Theta_\phi(t) \in \mathbb{R}$ admits a unique continuous extension $\bar{L}_t : L^2(\mathbb{T}^d) \ni \phi \mapsto \bar{\Theta}_\phi(t) \in \mathbb{R}$. Hence by Riesz representation theorem, there exists $\bar{\rho}^\nu(t,\cdot) \in L^2(\mathbb{T}^d)$ such that

$$\int_{\mathbb{T}^d} \bar{\rho}^\nu(t,x)\phi(x)dx = \bar{\Theta}_\phi(t) \quad \forall \phi \in L^2(\mathbb{T}^d), \quad \forall t \in [0, T],$$

with the properties that $\bar{\rho}^\nu(t,\cdot) = \rho^\nu(t,\cdot)$ for every $t \in [0, T] \setminus \mathcal{T}$, and $\|\bar{\rho}^\nu\|_{L^\infty((0,T);L^2(\mathbb{T}^d))} \leq \|\rho_{in}\|_{L^2(\mathbb{T}^d)}$.

We are left to prove that $[0, T] \ni t \longmapsto \Theta_\phi(t) = \int_{\mathbb{T}^d} \bar{\rho}^\nu(t, x)\phi(x)dx$ is continuous for any $\phi \in L^2(\mathbb{T}^d)$. By density of $\mathcal{N}$ in $L^2(\mathbb{T}^d)$, we that for every $\delta > 0$, there exists $\phi_\delta \in \mathcal{N}$ such that $\|\phi_\delta - \phi\|_{L^2(\mathbb{T}^d)} \leq \delta$. Hence, for every $t \in [0, T]$, we have

$$\left| \int_{\mathbb{T}^d} \bar{\rho}^\nu(t, x)\phi_\delta(x)dx - \int_{\mathbb{T}^d} \bar{\rho}^\nu(t, x)\phi(x)dx \right| \leq \|\rho^\nu\|_{L^\infty((0,T);L^2(\mathbb{T}^d))} \|\phi_\delta - \phi\|_{L^2(\mathbb{T}^d)}$$
$$\leq C_{T,\nu} \|\rho_{in}\|_{L^2(\mathbb{T}^d)} \delta,$$

where the constant $C_{T,\nu}$ depends on $T$ and $\nu$ only, and we have used (2.7) in the last inequality. Hence $[0, T] \ni t \longmapsto \int_{\mathbb{T}^d} \bar{\rho}^\nu(t, x)\phi_\delta(x)dx$ converges uniformly to $[0, T] \ni t \longmapsto \int_{\mathbb{T}^d} \bar{\rho}^\nu(t, x)\phi(x)dx$, which is continuous as well. This proves that $\bar{\rho}^\nu$ belongs to $C([0, T]; w - L^2(\mathbb{T}^d))$. As $\rho^\nu$ coincides with $\bar{\rho}^\nu$ except on a zero-measure set, we now assume that $\rho^\nu$ belongs to $C([0, T]; w - L^2(\mathbb{T}^d))$.

Mollifying ($\nu$ -PDE) by $w^\delta$, we then have

$$\begin{cases} \partial_t \rho^{\nu,\delta} - \nu \Delta \rho^{\nu,\delta} + \text{div}(\boldsymbol{b}\rho^{\nu,\delta}) = r^\delta, \\ \rho^{\nu,\delta}(0, \cdot) = \rho_{in}^\delta(\cdot), \end{cases} \quad (2.13)$$

with $r^\delta := \text{div}(\boldsymbol{b}\rho^{\nu,\delta} - (\boldsymbol{b}\rho^\nu) * w^\delta)$. Since for any $t \in [0, T]$ we have $\rho^\nu(t, \cdot) \in L^2(\mathbb{T}^d)$, it holds that

$$\|\rho^{\nu,\delta}(t, \cdot) - \rho^\nu(t, \cdot)\|_{L^2(\mathbb{T}^d)} \to 0$$

as $\delta \downarrow 0$. We can multiply the equation (2.13) by $\rho^{\nu,\delta}$ and integrate in space and in time over $(s, t) \times \mathbb{T}^d$ to get

$$\left| \|\rho^{\nu,\delta}(s, \cdot)\|^2_{L^2(\mathbb{T}^d)} - \|\rho^{\nu,\delta}(t, \cdot)\|^2_{L^2(\mathbb{T}^d)} \right| = 2 \left| \nu \int_s^t \|\nabla \rho^{\nu,\delta}(\tau, \cdot)\|^2_{L^2(\mathbb{T}^d)} d\tau + \int_s^t \int_{\mathbb{T}^d} r^\delta \rho^{\nu,\delta} \right|$$
$$\leq 2\nu \int_s^t \|\nabla \rho^\nu(\tau, \cdot)\|^2_{L^2(\mathbb{T}^d)} d\tau$$
$$+ 2\|r_\delta\|_{L^2((0,T);\dot{H}^{-1}(\mathbb{T}^d))} \|\rho^{\nu,\delta}\|_{L^2((0,T);\dot{H}^1(\mathbb{T}^d))}$$

We observe that $2\nu \|\rho^{\nu,\delta}\|^2_{L^2((0,T);\dot{H}^1(\mathbb{T}^d))} \leq \|\rho_{in}^\delta\|^2_{L^2(\mathbb{T}^d)} \leq \|\rho_{in}\|^2_{L^2(\mathbb{T}^d)}$ is uniformly bounded in $\delta > 0$ by (2.3) and that $\|r^\delta\|_{L^2((0,T);\dot{H}^{-1}(\mathbb{T}^d))} \to 0$ as $\delta \downarrow 0$ as noticed in the previous step. Then, passing to the limit $\delta \downarrow 0$ in the previous inequality, we get for any $0 \leq s \leq t \leq T$

$$\left| \|\rho^\nu(s, \cdot)\|^2_{L^2(\mathbb{T}^d)} - \|\rho^\nu(t, \cdot)\|^2_{L^2(\mathbb{T}^d)} \right| \leq 2\nu \int_s^t \|\nabla \rho^\nu(\tau, \cdot)\|^2_{L^2(\mathbb{T}^d)} d\tau.$$

Since $\rho^\nu \in L^2((0, T); H^1(\mathbb{T}^d))$, this inequality implies that $[0, T] \ni t \longmapsto \int_{\mathbb{T}^d} |\rho^\nu(t, x)|^2 dx$ is a continuous function. Since we also have $\rho^\nu \in C([0, T]; w - L^2(\mathbb{T}^d))$, we deduce $\rho^\nu \in C([0, T]; L^2(\mathbb{T}^d))$.

□

We record the following proposition, which will be used to prove the main theorem.

**Proposition 2.5** *Consider a divergence-free vector field $\boldsymbol{b} \in L^2((0, T) \times \mathbb{T}^d; \mathbb{R}^d)$, and an initial datum $\rho_{in} \in L^\infty(\mathbb{T}^d)$. Then there exists a vanishing diffusivity solution $\rho$ with initial datum $\rho_{in}$. Furthermore, any subsequence of $\{\rho^\nu\}_{\nu \in (0,1)}$ converging to $\rho$ weakly-star in $L^\infty((0, T) \times \mathbb{T}^d)$ is also converging in $C([0, T]; w - L^2(\mathbb{T}^d))$.*

***Proof* Step 1** (Existence) The family of unique weak solutions $\{\rho^\nu\}_{\nu>0} \subset L^\infty((0,T) \times \mathbb{T}^d)$ is relatively compact with respect to the weakly-star convergence, thanks to the estimate $\|\rho^\nu\|_{L^\infty((0,T)\times\mathbb{T}^d)} \leq \|\rho_{in}\|_{L^\infty(\mathbb{T}^d)}$ given in Theorem 2.4. Then, $\{\rho^\nu\}_{\nu>0}$ admits an accumulation point $\rho$, which solves (PDE) along $\boldsymbol{b}$ with initial datum $\rho_{in}$ and which is thus a vanishing diffusivity solution.

**Step 2** (Convergence in $C([0,T]; w - L^2(\mathbb{T}^d))$) Let $(\rho^{\nu_i})_{i\in\mathbb{N}}$ be a sequence, which converges to $\rho$ weak-star in $L^\infty((0,T) \times \mathbb{T}^d)$. Let us prove the convergence in $C([0,T]; w - L^2(\mathbb{T}^d))$. From Theorem 2.4, we have $\rho^\nu \in C([0,T]; L^2(\mathbb{T}^d))$. So, for every $\phi \in L^2(\mathbb{T}^d)$, we define the functions

$$f_\phi^\nu : [0,T] \ni t \longmapsto \int_{\mathbb{T}^d} \rho^\nu(t,x)\phi(x)dx.$$

Let $\mathcal{N} \subset C_c^2(\mathbb{T}^d)$ be a countable, dense set in $L^2(\mathbb{T}^d)$, and let $\phi \in \mathcal{N}$. We can prove that the family $\{f_\phi^\nu\}_{\nu>0}$ is bounded; thanks to the bound in Theorem 2.4, we have equiboundedness:

$$\left| \int_{\mathbb{T}^d} \rho^\nu(t,x)\phi(x)dx \right| \leq \|\rho_{in}\|_{L^\infty(\mathbb{T}^d)} \|\phi\|_{L^2(\mathbb{T}^d)}.$$

Also, we have the equicontinuity. Indeed, for every $0 \leq s \leq t \leq T$, by the weak formulation of ($\nu$ -PDE), we have

$$\left| f_\phi^\nu(t) - f_\phi^\nu(s) \right| \leq \int_s^t \int_{\mathbb{T}^d} |\boldsymbol{b}(\tau,x)\rho^\nu(\tau,x) \cdot \nabla\phi(x) + \nu\rho^\nu(\tau,x)\Delta\phi(x)|dxd\tau$$
$$\leq \left( \|\boldsymbol{b}\|_{L^1((s,t);L^2(\mathbb{T}^d))} \|\phi\|_{C^1(\mathbb{T}^d)} + \nu(t-s)\|\phi\|_{C^2(\mathbb{T}^d)} \right) \|\rho_{in}\|_{L^\infty(\mathbb{T}^d)}$$

where we have taken a test function $\phi(x)\psi(\tau) = \phi(x)\mathbb{1}_{[s,t]}(\tau)$, and the bound in Theorem 2.4. Therefore, by the Arzelà-Ascoli Theorem, and up to a relabeling of sequences, for every $\phi \in \mathcal{N}$, the sequence $(f_\phi^{\nu_i})_{i\in\mathbb{N}}$ converges in $C([0,T])$ to some $f_\phi$. We now claim that

$$f_\phi(t) = \int_{\mathbb{T}^d} \rho(t,x)\phi(x)dx \quad \forall t \in [0,T], \quad \forall \phi \in \mathcal{N} \tag{2.14}$$

where we recall that $\rho$ is the vanishing diffusivity solution. Let $\phi \in \mathcal{N}$, and observe that

$$\int_0^T f_\phi(s)\psi(s)ds = \int_0^T \int_{\mathbb{T}^d} \rho(s,x)\psi(s)\phi(x)dxds \quad \forall \psi \in L^\infty((0,T)) \tag{2.15}$$

since $\rho^{\nu_i}$ converges weakly-star to $\rho$ in $L^\infty((0,T) \times \mathbb{T}^d)$ as $i \to +\infty$. Therefore

$$f_\phi(s) = \int_{\mathbb{T}^d} \rho(s,x)\phi(x)dx \quad \text{for a.e. } s \in (0,T). \tag{2.16}$$

Since $f_\phi \in C([0,T])$, and by Theorem 2.4, $[0,T] \ni t \longmapsto \int_{\mathbb{T}^d} \rho(t,x)\phi(x)dx$ is continuous, the claim (2.14) holds $\forall t \in [0,T]$. Now, by density of $\mathcal{N}$ in $L^2(\mathbb{T}^d)$, we have that $\rho^{\nu_i}$ converges in $C([0,T]; w - L^2(\mathbb{T}^d))$ to $\rho$ as $i \to +\infty$. □

## 3 Analysis of backward problems

Throughout this section, we consider $\chi \in C_c^\infty((0,T) \times \mathbb{T}^d)$, a divergence-free vector field $\boldsymbol{b} : [0,T] \times \mathbb{T}^d \to \mathbb{R}^d$, which belongs to $L^1_{\text{loc}}((0,T]); BV(\mathbb{T}^d; \mathbb{R}^d)) \cap L^2((0,T) \times \mathbb{T}^d; \mathbb{R}^d))$, and $\nu \in (0,1)$. We consider the following backward advection problem

$$\begin{cases} \partial_t \theta_\chi + \text{div}(\boldsymbol{b}\theta_\chi) + \chi = 0, \\ \theta_\chi(T, \cdot) = 0. \end{cases} \tag{BW}$$

Weak solutions to (BW) are defined as follows.

**Definition 3.1** We shall say that $\theta_\chi \in L^\infty((0, T) \times \mathbb{T}^d)$ is a bounded weak solution of (BW), if for every $\phi \in C_c^\infty((0, T] \times \mathbb{T}^d)$, we have

$$\int_0^T \int_{\mathbb{T}^d} \Big[ \theta_\chi(t,x)\partial_t \phi(t,x) + \boldsymbol{b}(t,x)\theta_\chi(t,x) \cdot \nabla \phi(t,x) - \chi(t,x)\phi(t,x) \Big] dx dt = 0 \tag{3.1}$$

We consider also the following advection–diffusion problem.

$$\begin{cases} \partial_t \theta_\chi^\nu + \text{div}(\boldsymbol{b}\theta_\chi^\nu) + \nu \Delta \theta_\chi^\nu + \chi = 0, \\ \theta_\chi^\nu(T, \cdot) = 0. \end{cases} \tag{$\nu$–BW}$$

Similarly, weak solutions to ($\nu$–BW) are defined as follows.

**Definition 3.2** We shall say that $\theta_\chi^\nu \in L^\infty((0, T) \times \mathbb{T}^d)$ is a bounded weak solution of ($\nu$–BW), if for every $\phi \in C_c^\infty((0, T] \times \mathbb{T}^d)$, we have

$$\int_0^T \int_{\mathbb{T}^d} \Big[ \theta_\chi^\nu(t,x)\partial_t \phi(t,x) + \boldsymbol{b}(t,x)\theta_\chi^\nu(t,x) \cdot \nabla \phi(t,x) \\ - \nu \theta_\chi^\nu(t,x)\Delta \phi(t,x) - \chi(t,x)\phi(t,x) \Big] dx dt = 0. \tag{3.2}$$

### 3.1 The non-diffusive problem

The following result is a consequence of Theorem 2.1.

**Theorem 3.3** *In the context of this section, there exists a unique bounded weak solution of (BW) such that*

$$\|\theta_\chi\|_{L^\infty((0,T)\times\mathbb{T}^d)} \leq \int_0^T \|\chi(s,\cdot)\|_{L^\infty(\mathbb{T}^d)} ds \tag{3.3}$$

*and, up to redefining $\theta_\chi$ on a zero-measure set, we have $\theta_\chi \in C([0, T]; w - L^2(\mathbb{T}^d))$.*

*Remark 3.4* For the vector field constructed by Depauw in [13], which satisfies the assumptions of the above theorem, it is false that $\theta_\chi$ is continuous in time up to $t = 0$ in the strong $L^2(\mathbb{T}^d)$ topology, i.e. there exist $\chi$ and $\theta_\chi$ as in the above theorem such that $\theta_\chi$ does not belong to $C([0, T]; L^2(\mathbb{T}^d))$.

**Proof** **Step 1** (Existence) Let $\{w^\delta\}_{\delta>0}$ be a standard family of mollifiers; define $\boldsymbol{b}^\delta := \boldsymbol{b} * w^\delta$, which is such that $\boldsymbol{b}^\delta \to \boldsymbol{b}$ in $L^1((0, T) \times \mathbb{T}^d; \mathbb{R}^d)$, and consider the classical solution $\theta_\chi^\delta$ to (BW) with mollified vector field:

$$\begin{cases} \partial_t \theta_\chi^\delta + \text{div}(\boldsymbol{b}^\delta \theta_\chi^\delta) + \chi = 0, \\ \theta_\chi^\delta(T, \cdot) = 0. \end{cases} \tag{3.4}$$

The following bound holds true

$$\|\theta_\chi^\delta\|_{L^\infty((0,T)\times\mathbb{T}^d)} \leq \int_0^T \|\chi(s,\cdot)\|_{L^\infty(\mathbb{T}^d)} ds. \tag{3.5}$$

Existence of a solution $\theta_\chi$ to (BW) follows by the fact that the family $\{\theta_\chi^\delta\}_{\delta>0}$ is uniformly bounded in $L^\infty((0,T)\times\mathbb{T}^d)$, hence it admits an accumulation point in the weak-star topology. Using that $\boldsymbol{b}^\delta \to \boldsymbol{b}$ in $L^2((0,T)\times\mathbb{T}^d;\mathbb{R}^d)$, we conclude existence of a bounded weak solution $\theta_\chi$, since we can pass into the limit $\delta \downarrow 0$ in the weak formulation of ($\nu$–BW). Finally, (3.3) holds for $\theta_\chi$ thanks to the weak lower semicontinuity of the norm.

**Step 2** (Uniqueness). Recall that the equation is linear, hence the problem for the difference of two distinct solutions $v \in L^\infty((0,T)\times\mathbb{T}^d)$ reads as

$$\begin{cases} \partial_t v + \operatorname{div}(\boldsymbol{b} v) = 0, \\ v(T,\cdot) = 0. \end{cases}$$

Switching to the time variable $\tilde{t} = T-t$, we can refer to Theorem 2.1 to conclude $v=0$ and hence uniqueness.

**Step 3** (Time continuity) Let $\mathcal{N}$ be a dense, countable subset of $L^2(\mathbb{T}^d)$ consisting of smooth functions, which exists by separability, and let $\phi \in \mathcal{N}$. Let also $\psi \in C_c^\infty((0,T])$. Define

$$\Theta_\phi(t) := \int_{\mathbb{T}^d} \theta_\chi(t,x)\phi(x)dx,$$

and note that from the weak formulation of (BW):

$$\int_0^T \Theta_\phi(s)\partial_s\psi(s)ds = -\int_0^T \int_{\mathbb{T}^d} \boldsymbol{b}(\tau,x)\theta_\chi(\tau,x)\cdot\nabla\phi(x)\psi(\tau)dxd\tau$$
$$+ \int_0^T \int_{\mathbb{T}^d} \chi(\tau,x)\phi(x)\psi(\tau)dxd\tau.$$

Since $\boldsymbol{b} \in L^2((0,T)\times\mathbb{T}^d;\mathbb{R}^d)$, $\theta_\chi \in L^\infty((0,T);L^2(\mathbb{T}^d))$, $\nabla\phi \in L^\infty(\mathbb{T}^d)$ and $\psi \in L^\infty((0,T))$, we can conclude that $\partial_t\Theta_\phi \in L^1((0,T))$, hence $\Theta_\phi$ is continuous up to a set $\mathcal{T}_\phi$ of measure zero. Now let

$$\mathcal{T} := \bigcup_{\phi\in\mathcal{N}} \mathcal{T}_\phi,$$

which also has zero measure. For every $\phi \in \mathcal{N}$, we redefine $\Theta_\phi$ on $\mathcal{T}$ to be continuous, and we observe that, since $\Theta_\phi(t) \leq \|\theta_\chi\|_{L^\infty((0,T);L^2(\mathbb{T}^d))}\|\phi\|_{L^2(\mathbb{T}^d)}$ for all $t \in [0,T]$, the bounded linear functional $L_t : \mathcal{N} \ni \phi \longmapsto \Theta_\phi(t) \in \mathbb{R}$ admits a unique continuous extension $\bar{L}_t : L^2(\mathbb{T}^d) \ni \phi \longmapsto \bar{\Theta}_\phi(t) \in \mathbb{R}$. Hence, by Riesz representation theorem, there exists $\bar{\theta}_\chi(t,\cdot) \in L^2(\mathbb{T}^d)$ such that

$$\int_{\mathbb{T}^d} \bar{\theta}_\chi(t,x)\phi(x)dx = \bar{\Theta}_\phi(t) \quad \forall \phi \in L^2(\mathbb{T}^d), \ \forall t \in [0,T],$$

with the properties $\bar{\theta}_\chi(t,\cdot) = \theta_\chi(t,\cdot) \ \forall t \in [0,T] \setminus \mathcal{T}$ and $\|\bar{\theta}_\chi\|_{L^\infty((0,T);L^2(\mathbb{T}^d))} \leq \|\theta_\chi\|_{L^\infty((0,T);L^2(\mathbb{T}^d))}$.

We are left to prove that $[0,T] \ni t \longmapsto \bar{\Theta}_\phi(t) = \int_{\mathbb{T}^d} \bar{\theta}_\chi(t,x)\phi(x)dx$ is continuous for any $\phi \in L^2(\mathbb{T}^d)$. By density of $\mathcal{N}$ in $L^2(\mathbb{T}^d)$, we know that $\forall \delta > 0, \exists \phi_\delta \in \mathcal{N}$ such that $\|\phi_\delta - \phi\|_{L^2(\mathbb{T}^d)} \leq \delta$. Hence, for every $t \in [0,T]$, we have

$$\left|\int_{\mathbb{T}^d} \bar{\theta}_\chi(t,x)\phi_\delta(x)dx - \int_{\mathbb{T}^d} \bar{\theta}_\chi(t,x)\phi(x)dx\right| \leq \|\theta_\chi\|_{L^\infty((0,T);L^2(\mathbb{T}^d))}\|\phi_\delta - \phi\|_{L^2(\mathbb{T}^d)}$$
$$\leq \|\chi\|_{L^1((0,T);L^\infty(\mathbb{T}^d))}\delta,$$

where we have used (3.3) in the last inequality. Hence, $[0, T] \ni t \longmapsto \int_{\mathbb{T}^d} \bar{\theta}_\chi(t, x)\phi_\delta(x)dx$ converges uniformly to $[0, T] \ni t \longmapsto \int_{\mathbb{T}^d} \bar{\theta}_\chi(t, x)\phi(x)dx$ as $\delta \downarrow 0$, which implies that $[0, T] \ni t \longmapsto \int_{\mathbb{T}^d} \bar{\theta}_\chi(t, x)\phi(x)dx$ is continuous as well. This proves that $\bar{\theta}_\chi$ belongs to $C([0, T]; w - L^2(\mathbb{T}^d))$. As $\theta_\chi$ coincides with $\bar{\theta}_\chi$ except on a zero-measure set, the thesis follows.

□

### 3.2 The diffusive problem

Let us summarise the classical theory for ($\nu$−BW) along a smooth vector field.

**Remark 3.5** Let $(\Omega, (\mathcal{F}_t)_{t\geq 0}, \mathbb{P})$ be a filtered probability space. For $t \in [0, T]$, consider the following forward stochastic differential equation

$$dX^\nu_{s,t} = \boldsymbol{b}(s, X^\nu_{s,t})ds + \sqrt{2\nu}d\boldsymbol{W}_s \quad \text{with} \quad X^\nu_{t,t}(x) = x, \quad (3.6)$$

where $\boldsymbol{W}$ is a $\mathbb{T}^d$-valued Brownian motion adapted to the filtration such that we have $\boldsymbol{W}_t = 0$ almost surely. Assume now that the vector field $\boldsymbol{b}$ is divergence-free and smooth. Then for almost every $\omega \in \Omega$, we have $X_{s,t}(\omega)$ is incompressible, namely $X_{s,t}(\omega, \cdot)_\# \mathscr{L}^d = \mathscr{L}^d$. We also recall that the Feynman-Kac formula, which gives the following expression for the solution of ($\nu$−BW)

$$\theta^\nu_\chi(t, x) = \int_t^T \mathbb{E}[\chi(s, X^\nu_{s,t}(x))]ds.$$

It can be checked that $\theta^\nu_\chi \in C([0, T]; L^2(\mathbb{T}^d))$, and a computation using the dual characterisation of the norm and Fubini shows that

$$\sup_{t\in[0,T]} \|\theta^\nu_\chi(t, \cdot)\|_{L^\infty(\mathbb{T}^d)} \leq \int_0^T \|\chi(s, \cdot)\|_{L^\infty(\mathbb{T}^d)}ds.$$

The following lemma will be useful.

**Lemma 3.6** *In the context of this section, there exists a unique bounded weak solution $\theta^\nu_\chi \in C([0, T]; L^2(\mathbb{T}^d)) \cap L^2((0, T); H^1(\mathbb{T}^d))$ of ($\nu$−BW), which satisfies the following energy estimate for every $t \in [0, T]$*

$$\|\theta^\nu_\chi(t, \cdot)\|^2_{L^2(\mathbb{T}^d)} + 2\nu \int_t^T \int_{\mathbb{T}^d} |\nabla \theta^\nu_\chi(s, x)|^2 dxds \leq 4T \|\chi\|^2_{L^2((t,T)\times\mathbb{T}^d)}. \quad (3.7)$$

*and the bound*

$$\|\theta^\nu_\chi\|_{L^\infty((0,T)\times\mathbb{T}^d)} \leq \int_0^T \|\chi(s, \cdot)\|_{L^\infty(\mathbb{T}^d)}ds. \quad (3.8)$$

**Proof** **Step 1** (Existence) Consider a standard family of mollifiers $\{w^\delta\}_{\delta>0}$ and define $\boldsymbol{b}^\delta = \boldsymbol{b} * w^\delta$. Define $\theta^{\nu,\delta}_\chi$ as the solution to the backwards problem with mollified data:

$$\begin{cases} \partial_t \theta^{\nu,\delta}_\chi + \nu\Delta\theta^{\nu,\delta}_\chi + \text{div}(\boldsymbol{b}^\delta \theta^{\nu,\delta}_\chi) + \chi = 0, \\ \theta^{\nu,\delta}_\chi(T, \cdot) = 0. \end{cases} \quad (3.9)$$

By the classical theory detailed in Remark 3.5, there exists a unique smooth solution $\theta^{\nu,\delta}$; now, mutliply (3.9) by $\theta^{\nu,\delta}_\chi$ and integrate in space and time over $(t, T) \times \mathbb{T}^d$ to get for every $t \in [0, T]$ that

$$\|\theta_\chi^{\nu,\delta}(t,\cdot)\|_{L^2(\mathbb{T}^d)}^2 + 2\nu \int_t^T \int_{\mathbb{T}^d} |\nabla \theta_\chi^{\nu,\delta}(s,x)|^2 dx ds \leq 4T \|\chi\|_{L^2((0,T)\times\mathbb{T}^d)}^2. \quad (3.10)$$

We also have that

$$\|\theta_\chi^{\nu,\delta}\|_{L^\infty((0,T)\times\mathbb{T}^d)} \leq \int_0^T \|\chi(s,\cdot)\|_{L^\infty(\mathbb{T}^d)} ds. \quad (3.11)$$

Then, by compactness arguments, there exists a subsequence $(\theta_\chi^{\nu,\delta_k})_{k\in\mathbb{N}}$ such that $\theta_\chi^{\nu,\delta_k}$ converges weakly-star in $L^\infty((0,T)\times\mathbb{T}^d)$ to $\theta_\chi^\nu$, and converges weakly in $L^2((0,T); H^1(\mathbb{T}^d))$. Moreover $\theta_\chi^\nu$ is a bounded weak solution of ($\nu$−BW), since the equation is linear. The bounds (3.7) and (3.8) for $\theta_\chi^\nu$ follow from the weak lower semicontinuity of the norms.

**Step 2** (Uniqueness) We observe that the difference of two solutions $v \in L^2((0,T); H^1(\mathbb{T}^d)) \cap C([0,T]; L^2(\mathbb{T}^d))$ solves

$$\begin{cases} \partial_t v + \nu \Delta v + \mathrm{div}(\boldsymbol{b} v) = 0, \\ v(T,\cdot) = 0. \end{cases}$$

We can conclude uniqueness by mollification, analogously to the case of the forward problem (2.8). Indeed, introducing the mollified vector field $\boldsymbol{b}^\delta = \boldsymbol{b} * w^\delta$, we have that the mollified solutions $v^\delta = v * w^\delta$ solves

$$\begin{cases} \partial_t v^\delta + \nu \Delta v^\delta + \mathrm{div}(\boldsymbol{b}^\delta v^\delta) = r^\delta, \\ v(T,\cdot) = 0. \end{cases}$$

where $r^\delta := \mathrm{div}(\boldsymbol{b} v^\delta - (\boldsymbol{b} v) * w^\delta)$. Let $t \in [0,T]$. Now, multiplying by $v^\delta$ and integrating in space and time, we get

$$\frac{1}{2} \int_t^T \partial_s \int_{\mathbb{T}^d} v^\delta(s,x) v^\delta(s,x) dx ds - \nu \int_t^T \int_{\mathbb{T}^d} |\nabla v^\delta(s,x)|^2 dx ds$$
$$= \int_t^T \int_{\mathbb{T}^d} r^\delta(s,x) v^\delta(s,x) dx ds,$$

We then have

$$\frac{\|v^\delta(t,\cdot)\|_{L^2(\mathbb{T}^d)}^2}{2} + \nu \|v^\delta\|_{L^2((t,T);\dot{H}^1(\mathbb{T}^d))} \leq \|r^\delta\|_{L^2((t,T);\dot{H}^{-1}(\mathbb{T}^d))} \|v^\delta\|_{L^2((t,T);\dot{H}^1(\mathbb{T}^d))}$$

which implies, using Young inequality as in (2.12) that

$$\frac{\|v^\delta(t,\cdot)\|_{L^2(\mathbb{T}^d)}^2}{2} \leq \frac{1}{4\nu} \|r^\delta\|_{L^2((0,T);\dot{H}^{-1}(\mathbb{T}^d))}^2 \to 0 \quad \text{as } \delta \downarrow 0,$$

as done in Theorem 2.4, from which we conclude $v = 0$ a.e. because $t$ was arbitrary in $[0,T]$.

**Step 3** (Time continuity) The time continuity can be proved analogously to **Step 3** of Theorem 2.4.

□

**Lemma 3.7** *In the context of this section, consider $\rho^\nu \in C([0,T]; L^2(\mathbb{T}^d)) \cap L^2((0,T); H^1(\mathbb{T}^d))$ the unique weak solution to ($\nu$ -PDE), and $\theta_\chi^\nu \in C([0,T]; L^2(\mathbb{T}^d)) \cap L^2((0,T); H^1(\mathbb{T}^d))$ the unique bounded weak solution to ($\nu$−BW). Then the following duality formula holds*

$$\int_0^T \int_{\mathbb{T}^d} \rho^\nu(s,x) \chi(s,x) dx ds = \int_{\mathbb{T}^d} \rho_{in}(x) \theta_\chi^\nu(0,x) dx. \quad (3.12)$$

**Proof** We mollify the vector field and the initial datum in ($\nu$ -PDE), and consider $\rho^{\nu,\delta}$ solving

$$\begin{cases} \partial_t \rho^{\nu,\delta} + \text{div}(\boldsymbol{b}^\delta \rho^{\nu,\delta}) - \nu \Delta \rho^{\nu,\delta} = 0, \\ \rho^{\nu,\delta}(0, \cdot) = \rho_{in}^\delta. \end{cases} \tag{3.13}$$

Consider also:

$$\begin{cases} \partial_t \theta_\chi^\nu + \text{div}(\boldsymbol{b}\theta_\chi^\nu) + \nu \Delta \theta_\chi^\nu + \chi = 0, \\ \theta_\chi^\nu(T, \cdot) = 0. \end{cases} \tag{3.14}$$

Since $\rho^{\nu,\delta}$ is smooth, we can multiply the equation in (3.14) by $\rho^{\nu,\delta}$ and integrate by parts discharging derivatives on $\rho^{\nu,\delta}$:

$$- \int_0^T \int_{\mathbb{T}^d} (\partial_t \rho^{\nu,\delta}(t,x) + \boldsymbol{b}^\delta(t,x) \cdot \nabla \rho^{\nu,\delta}(t,x) - \nu \Delta \rho^{\nu,\delta}(t,x))\theta_\chi^\nu(t,x) dt dx$$
$$- \int_0^T \int_{\mathbb{T}^d} (\boldsymbol{b}(t,x) - \boldsymbol{b}^\delta(t,x)) \cdot \nabla \rho^{\nu,\delta}(t,x)\theta_\chi^\nu(t,x) dx dt$$
$$+ \int_0^T \int_{\mathbb{T}^d} \chi(t,x) \rho^{\nu,\delta}(t,x) dx dt$$
$$= \int_{\mathbb{T}^d} \rho_{in}^\delta(x) \theta_\chi^\nu(0,x) dx.$$

By (3.13) we have that $\partial_t \rho^{\nu,\delta} + \boldsymbol{b}^\delta \cdot \nabla \rho^{\nu,\delta} - \nu \Delta \rho^{\nu,\delta} = 0$; moreover, we can prove that, as $\delta \downarrow 0$:

$$\int_0^T \int_{\mathbb{T}^d} (\boldsymbol{b}(t,x) - \boldsymbol{b}^\delta(t,x)) \cdot \nabla \rho^{\nu,\delta}(t,x) \theta_\chi^\nu(t,x) dx dt$$
$$\leq \|\boldsymbol{b} - \boldsymbol{b}^\delta\|_{L^2((0,T)\times\mathbb{T}^d)} \|\nabla \rho^{\nu,\delta}\|_{L^2((0,T)\times\mathbb{T}^d)} \|\theta_\chi^\nu\|_{L^\infty((0,T)\times\mathbb{T}^d)} \to 0,$$

thanks to the fact that $\boldsymbol{b}^\delta \to \boldsymbol{b}$ strongly in $L^2((0,T) \times \mathbb{T}^d; \mathbb{R}^d)$, and by Lemma 3.6, $\|\theta_\chi^\nu\|_{L^\infty((0,T)\times\mathbb{T}^d)}$ is bounded. Note also that

$$\int_0^T \int_{\mathbb{T}^d} \chi(t,x) \rho^{\nu,\delta}(t,x) dx dt \to \int_0^T \int_{\mathbb{T}^d} \chi(t,x) \rho^\nu(t,x) dx dt \quad \text{as } \delta \downarrow 0,$$

since $\rho^{\nu,\delta} \overset{*}{\rightharpoonup} \rho^\nu$ in $L^\infty((0,T) \times \mathbb{T}^d)$ as proved in Theorem 2.4, and

$$\int_{\mathbb{T}^d} \rho_{in}^\delta(x) \theta_\chi^\nu(0,x) dx \to \int_{\mathbb{T}^d} \rho_{in}(x) \theta_\chi^\nu(0,x) dx \quad \text{as } \delta \downarrow 0,$$

because $\rho_{in}^\delta \to \rho_{in}$ in $L^2(\mathbb{T}^d)$. Therefore, we have that

$$\int_0^T \int_{\mathbb{T}^d} \chi(t,x) \rho^\nu(t,x) dx dt = \int_{\mathbb{T}^d} \rho_{in}(x) \theta_\chi^\nu(0,x) dx.$$

$\square$

### 3.3 Uniqueness of the backward vanishing diffusivity solution

The following lemma is essential to prove Theorem 1.4.

**Lemma 3.8** *In the context of this section, consider the family $\{\theta_\chi^\nu\}_{\nu>0}$ of bounded weak solutions of ($\nu$-BW), and the bounded weak solution $\theta_\chi$ of (BW). Then the family $\{\theta_\chi^\nu\}_{\nu>0}$ converges in $C([0, T]; w - L^2(\mathbb{T}^d))$ to $\theta_\chi$ as $\nu \downarrow 0$.*

*Proof* Since from Lemma 3.6, we have the uniform bound

$$\|\theta_\chi^\nu\|_{L^\infty((0,T)\times\mathbb{T}^d)} \leq \int_0^T \|\chi(s,\cdot)\|_{L^\infty(\mathbb{T}^d)} ds, \tag{3.15}$$

the family $\{\theta_\chi^\nu\}_{\nu>0}$ is bounded and therefore has at least one accumulation point in the weak-star topology on $L^\infty((0, T) \times \mathbb{T}^d)$. Moreover, any limit point must be a weak solution of (BW), and as by Theorem 3.3, $\theta_\chi$ is the unique such solution, we have that $\theta_\chi^\nu$ converges weakly-star in $L^\infty((0, T) \times \mathbb{T}^d)$ to $\theta_\chi$ as $\nu \downarrow 0$.

In particular, for every $\phi \in L^2(\mathbb{T}^d)$, and every $\psi \in L^1((0, T))$ it holds:

$$\int_0^T \int_{\mathbb{T}^d} \theta_\chi^\nu(s,x)\phi(x)\psi(s) dx ds \to \int_0^T \int_{\mathbb{T}^d} \theta_\chi(s,x)\phi(x)\psi(s) dx ds \quad \text{as } \nu \downarrow 0. \tag{3.16}$$

Let us prove that the convergence takes place in $C([0, T]; w - L^2(\mathbb{T}^d))$. From Lemma 3.6, we have $\theta_\chi^\nu \in C([0, T]; L^2(\mathbb{T}^d))$. So, for every $\phi \in L^2(\mathbb{T}^d)$, we define the functions

$$f_\phi^\nu : [0, T] \ni t \longmapsto \int_{\mathbb{T}^d} \theta_\chi^\nu(t,x)\phi(x) dx.$$

Let $\mathcal{N} \subset C_c^2(\mathbb{T}^d)$ be a countable, dense set in $L^2(\mathbb{T}^d)$, and let $\phi \in \mathcal{N}$. We can prove that the sequence $\{f_\phi^\nu\}_{\nu>0}$ is bounded; thanks to (3.7), we have equiboundedness:

$$\left| \int_{\mathbb{T}^d} \theta_\chi^\nu(t,x)\phi(x) dx \right| \leq \|\theta_\chi^\nu\|_{L^\infty((0,T);L^2(\mathbb{T}^d))} \|\phi\|_{L^\infty(\mathbb{T}^d)} \leq 4T \|\chi\|_{L^2((0,T)\times\mathbb{T}^d)} \|\phi\|_{L^\infty(\mathbb{T}^d)}.$$

Also, we have the equicontinuity. Indeed, for every $0 \leq s \leq t \leq T$ by the weak formulation of (3.2), we have

$$\left| f_\phi^\nu(t) - f_\phi^\nu(s) \right| \leq \int_s^t \int_{\mathbb{T}^d} |\boldsymbol{b}(\tau,x)\theta_\chi^\nu(\tau,x) \cdot \nabla\phi(x) + \nu\theta_\chi^\nu(\tau,x)\Delta\phi(x) + \chi(\tau,x)\phi(x)| dx d\tau$$

$$\leq \left( \|\boldsymbol{b}\|_{L^1((s,t);L^2(\mathbb{T}^d))} \|\phi\|_{C^1(\mathbb{T}^d)} + \nu(t-s)\|\phi\|_{C^2(\mathbb{T}^d)} \right) \|\theta_\chi^\nu\|_{L^\infty((0,T);L^2(\mathbb{T}^d))}$$

$$+ \|\chi\|_{L^1((s,t);L^1(\mathbb{T}^d))} \|\phi\|_{C(\mathbb{T}^d)}$$

$$\leq \left( \|\boldsymbol{b}\|_{L^1((s,t);L^2(\mathbb{T}^d))} \|\phi\|_{C^1(\mathbb{T}^d)} + \nu(t-s)\|\phi\|_{C^2(\mathbb{T}^d)} \right) 4T\|\chi\|_{L^2((0,T)\times\mathbb{T}^d)}$$

$$+ \|\chi\|_{L^1((s,t);L^1(\mathbb{T}^d))} \|\phi\|_{C(\mathbb{T}^d)},$$

where we have taken a test function $\phi(x)\psi(\tau) = \phi(x)1\!\!1_{[s,t]}(\tau)$, and used (3.15) in the last bound. Therefore, by the Arzelà-Ascoli theorem and a diagonal argument, there exists a sequence $(\nu_i)_{i\in\mathbb{N}}$ such that $\nu_i \downarrow 0$, and for every $\phi \in \mathcal{N}$, the sequence $(f_\phi^{\nu_i})_{i\in\mathbb{N}}$ converges in $C([0, T])$ to some $f_\phi$. Now, let $\phi \in \mathcal{N}$. We claim that

$$f_\phi(t) = \int_{\mathbb{T}^d} \theta_\chi(t,x)\phi(x) dx \quad \forall t \in [0, T], \tag{3.17}$$

where $\theta_\chi \in C([0, T]; w - L^2(\mathbb{T}^d))$ is the unique solution to the transport equation by Theorem 3.3. We observe by (3.23) that

$$\int_0^T f_\phi(s)\psi(s) ds = \int_0^T \int_{\mathbb{T}^d} \theta_\chi(s,x)\psi(s)\phi(x) dx ds \quad \forall \psi \in L^1((0, T)). \tag{3.18}$$

Therefore

$$f_\phi(s) = \int_{\mathbb{T}^d} \theta_\chi(s, x)\phi(x)dx \quad \text{for a.e. } s \in (0, T). \tag{3.19}$$

Since $f_\phi \in C([0, T])$, and by Lemma 3.6, $[0, T] \ni t \longmapsto \int_{\mathbb{T}^d} \theta_\chi(s, x)\phi(x)dx$ is continuous, the claim (3.24) holds $\forall t \in [0, T]$.

Therefore, the family $\{f_\phi^\nu\}_{\nu>0}$ converges uniformly to $f_\phi$, because every subsequence $\{f_\phi^{\nu_k}\}_{k \in \mathbb{N}}$ admits a subsequence $\{f_\phi^{\nu_{k_l}}\}_{l \in \mathbb{N}}$ convergent to the same limit $f_\phi$. Indeed, if we suppose by contradiction that the full sequence is not convergent, i.e., $\exists \varepsilon > 0, \{f_\phi^{\nu_k}\}_{k \in \mathbb{N}}$ such that $\|f_\phi^{\nu_k} - f_\phi\|_{C([0,T])} > \varepsilon$, we contradict the hypothesis. Finally, by density of $\mathcal{N}$ in $L^2(\mathbb{T}^d)$, we have that $\theta_\chi^\nu$ converges in $C([0, T]; w - L^2(\mathbb{T}^d))$ to $\theta_\chi$ as $\nu \downarrow 0$. □

### 3.4 Regularisation by convolution

Consider a standard mollifier $w \in C_c^\infty((0, 1)^d)$ such that $w \geq 0$ and $\int_{(0,1)^d} w(x)dx = 1$, and for every $\delta > 0$, write $w^\delta(x) = w(x/\delta)/\delta^d$, and $\boldsymbol{b}^\delta = \boldsymbol{b} * w^\delta$. Let $\theta_\chi^\delta$ be the unique solution

$$\begin{cases} \partial_t \theta_\chi^\delta + \operatorname{div}(\boldsymbol{b}^\delta \theta_\chi^\delta) + \chi = 0, \\ \theta_\chi^\delta(T, \cdot) = 0. \end{cases} \tag{$\delta$-BW}$$

It is then given by

$$\theta_\chi^\delta(t, x) = \int_t^T \chi(s, X^\delta(s, t, x))ds, \tag{3.20}$$

where $X^\delta$ solves

$$\begin{cases} \partial_s X^\delta(s, t, x) = \boldsymbol{b}^\delta(t, X^\delta(s, t, x)), \\ X^\delta(t, t, x) = x. \end{cases} \tag{$\delta$-ODE}$$

Let $\rho^\delta$ be the unique bounded weak solution of (PDE) along $\boldsymbol{b}^\delta$ with initial datum $\rho_{in}$. The following then holds.

**Lemma 3.9** *In the context of this section, consider the family $\{\theta_\chi^\delta\}_{\delta>0}$ of bounded weak solutions of ($\delta$-BW), and the bounded weak solution $\theta_\chi$ of (BW). Then the family $\{\theta_\chi^\delta\}_{\delta>0}$ converges in $C([0, T]; w - L^2(\mathbb{T}^d))$ to $\theta_\chi$ as $\delta \downarrow 0$.*

*Moreover, we have the duality formula*

$$\int_0^T \int_{\mathbb{T}^d} \rho^\delta(s, x)\chi(s, x)dxds = \int_{\mathbb{T}^d} \rho_{in}(x)\theta_\chi^\delta(0, x)dx \quad \forall \chi \in C_c^\infty((0, T) \times \mathbb{T}^d), \tag{3.21}$$

*for every $\delta > 0$.*

*Proof* The argument is analogous to the proof of Lemma 3.8. Thanks to the representation formula (3.20), we have the uniform in $\delta$ bound

$$\|\theta_\chi^\delta\|_{L^\infty((0,T)\times\mathbb{T}^d)} \leq \int_0^T \|\chi(s, \cdot)\|_{L^\infty(\mathbb{T}^d)}ds, \tag{3.22}$$

Thus, the family $\{\theta_\chi^\delta\}_{\delta>0}$ is bounded and therefore has at least one accumulation point for the weak-star topology on $L^\infty((0, T) \times \mathbb{T}^d))$. Moreover, any limit point must be a weak solution of (BW), and as by Theorem 3.3, $\theta_\chi$ is the unique such solution, we have that $\theta_\chi^\delta$ converges weakly-star in $L^\infty((0,T) \times \mathbb{T}^d)$ to $\theta_\chi$ as $\delta \downarrow 0$.

In particular, for every $\phi \in L^2(\mathbb{T}^d)$, and every $\psi \in L^1((0,T))$ it holds:

$$\int_0^T \int_{\mathbb{T}^d} \theta_\chi^\delta(s,x)\phi(x)\psi(s)dxds \to \int_0^T \int_{\mathbb{T}^d} \theta_\chi(s,x)\phi(x)\psi(s)dxds \quad \text{as } \delta \downarrow 0. \quad (3.23)$$

Let us prove that the convergence takes place in $C([0,T]; w - L^2(\mathbb{T}^d))$. From Lemma 3.6, we have $\theta_\chi^\delta \in C([0,T]; L^2(\mathbb{T}^d))$. So, for every $\phi \in L^2(\mathbb{T}^d)$, we define the functions

$$f_\phi^\delta : [0,T] \ni t \longmapsto \int_{\mathbb{T}^d} \theta_\chi^\delta(t,x)\phi(x)dx.$$

Let $\mathcal{N} \subset C_c^2(\mathbb{T}^d)$ be a countable, dense set in $L^2(\mathbb{T}^d)$, and let $\phi \in \mathcal{N}$. We can prove that the sequence $\{f_\phi^\delta\}_{\delta > 0}$ is bounded; thanks to (3.7), we have equiboundedness:

$$\left| \int_{\mathbb{T}^d} \theta_\chi^\delta(t,x)\phi(x)dx \right| \leq 4T \|\theta_\chi^\delta\|_{L^\infty((0,T);L^2(\mathbb{T}^d))} \|\phi\|_{L^\infty(\mathbb{T}^d)} \leq 4T \|\chi\|_{L^2((0,T)\times\mathbb{T}^d)} \|\phi\|_{L^\infty(\mathbb{T}^d)}.$$

Also, we have the equicontinuity. Indeed, for every $0 \leq s \leq t \leq T$ by the weak formulation of (3.2), we have

$$\left| f_\phi^\delta(t) - f_\phi^\delta(s) \right| \leq \int_s^t \int_{\mathbb{T}^d} |\boldsymbol{b}(\tau,x)\theta_\chi^\delta(\tau,x) \cdot \nabla\phi(x) + \chi(\tau,x)\phi(x)|dxd\tau$$

$$\leq \|\boldsymbol{b}\|_{L^1((s,t);L^2(\mathbb{T}^d))} \|\phi\|_{C^1(\mathbb{T}^d)} \|\theta_\chi^\delta\|_{L^\infty((0,T);L^2(\mathbb{T}^d))}$$
$$+ \|\chi\|_{L^1((s,t);L^1(\mathbb{T}^d))} \|\phi\|_{C(\mathbb{T}^d)}$$
$$\leq 4T \|\boldsymbol{b}\|_{L^1((s,t);L^2(\mathbb{T}^d))} \|\phi\|_{C^1(\mathbb{T}^d)} \|\chi\|_{L^2((0,T)\times\mathbb{T}^d)}$$
$$+ \|\chi\|_{L^1((s,t);L^1(\mathbb{T}^d))} \|\phi\|_{C(\mathbb{T}^d)},$$

where we have taken a test function $\phi(x)\psi(\tau) = \phi(x)\mathbb{1}_{[s,t]}(\tau)$, and used (3.15) in the last bound. Therefore, by the Arzelà-Ascoli theorem and a diagonal argument, there exists a sequence $(\delta_i)_{i \in \mathbb{N}}$ such that $\delta_i \downarrow 0$, and for every $\phi \in \mathcal{N}$, the sequence $(f_\phi^{\delta_i})_{i \in \mathbb{N}}$ converges in $C([0,T])$ to some $f_\phi$. Now, let $\phi \in \mathcal{N}$. We claim that

$$f_\phi(t) = \int_{\mathbb{T}^d} \theta_\chi(t,x)\phi(x)dx \quad \forall t \in [0,T], \quad (3.24)$$

where $\theta_\chi \in C([0,T]; w - L^2(\mathbb{T}^d))$ is the unique solution to the transport equation by Theorem 3.3. We observe by (3.23) that

$$\int_0^T f_\phi(s)\psi(s)ds = \int_0^T \int_{\mathbb{T}^d} \theta_\chi(s,x)\psi(s)\phi(x)dxds \quad \forall \psi \in L^1((0,T)). \quad (3.25)$$

Therefore

$$f_\phi(s) = \int_{\mathbb{T}^d} \theta_\chi(s,x)\phi(x)dx \quad \text{for a.e. } s \in (0,T). \quad (3.26)$$

Since $f_\phi \in C([0,T])$, and by Lemma 3.6, $[0,T] \ni t \longmapsto \int_{\mathbb{T}^d} \theta_\chi(s,x)\phi(x)dx$ is continuous, the claim (3.24) holds $\forall t \in [0,T]$.

Therefore, the family $\{f_\phi^\delta\}_{\delta > 0}$ converges uniformly to $f_\phi$, because every subsequence $\{f_\phi^{\delta_k}\}_{k \in \mathbb{N}}$ admits a subsequence $\{f_\phi^{\delta_{k_l}}\}_{l \in \mathbb{N}}$ convergent to the same limit $f_\phi$. Indeed, if we suppose by contradiction that the full sequence is not convergent, i.e., $\exists \varepsilon > 0$, $\{f_\phi^{\delta_k}\}_{k \in \mathbb{N}}$ such that $\|f_\phi^{\delta_k} - f_\phi\|_{C([0,T])} > \varepsilon$, we contradict the hypothesis. Finally, by density of $\mathcal{N}$ in $L^2(\mathbb{T}^d)$, we have that $\theta_\chi^\delta$ converges in $C([0,T]; w - L^2(\mathbb{T}^d))$ to $\theta_\chi$ as $\delta \downarrow 0$.

Let us prove the duality formula (3.21). Let $\chi \in C_c^\infty((0,T) \times \mathbb{T}^d)$ be arbitrary, and note that by the classical theory

$$\rho^\delta(s, \cdot) = X^\delta(s, 0, \cdot)_\# \rho_{in} \mathscr{L}^d, \qquad (3.27)$$

where recall that $X^\delta$ solves ($\delta$-ODE). Therefore,

$$\int_0^T \int_{\mathbb{T}^d} \rho^\delta(s,x)\chi(s,x)dxds = \int_0^T \int_{\mathbb{T}^d} \rho_{in}(x)\chi(s, X^\delta(s,0,x))dxds$$
$$= \int_{\mathbb{T}^d} \rho_{in}(x)\theta_\chi^\delta(0,x)dx, \qquad (3.28)$$

where in the last equality, we have used Fubini, and (3.20). This proves the thesis. □

## 4 Proof of Theorem 1.4

Let $\chi \in C_c^\infty((0,T) \times \mathbb{T}^d)$ be arbitrary. We have proved the duality formula in (3.12)

$$\int_0^T \int_{\mathbb{T}^d} \rho^\nu(s,x)\chi(s,x)dxds = \int_{\mathbb{T}^d} \rho_{in}(x)\theta_\chi^\nu(0,x)dx. \qquad (4.1)$$

Let $\rho$ be a vanishing viscosity solution of (PDE) along $\boldsymbol{b}$ with initial datum $\rho_{in}$, which exists by Proposition 2.5. We can pass into the limit $\nu \downarrow 0$ in (3.12) and we get

$$\int_0^T \int_{\mathbb{T}^d} \rho(s,x)\chi(s,x)dxds = \int_{\mathbb{T}^d} \rho_{in}(x)\theta_\chi(0,x)dx, \qquad (4.2)$$

where we have used that $\theta_\chi^\nu(0,\cdot) \rightharpoonup \theta_\chi(0,\cdot)$ in $L^2(\mathbb{T}^d)$ as $\nu \downarrow 0$, thanks to Lemma 3.8. As $\chi$ was arbitrary, this uniquely characterises the vanishing diffusivity solution, up to a zero measure set.

Now let $w$ be a standard mollifier, and let $\{\tilde{\rho}^\delta\}_{\delta > 0}$ be the family of unique bounded weak solutions along $\boldsymbol{b} * w^\delta$. Then by the classical theory $\|\tilde{\rho}^\delta\|_{L^\infty((0,T) \times \mathbb{T}^d)} \leq \|\rho_{in}\|_{L^\infty(\mathbb{T}^d)}$ for every $\delta > 0$, and there is therefore a sequence $(\delta_i)_{i \in \mathbb{N}}$ such that $\tilde{\rho}^{\delta_i}$ converges weakly-star in $L^\infty((0,T) \times \mathbb{T}^d)$ to some $\tilde{\rho}$ as $i \to +\infty$. Also by Lemma 3.9, we have the duality formula

$$\int_0^T \int_{\mathbb{T}^d} \tilde{\rho}^\delta(t,x)\chi(t,x)dxdt = \int_0^T \int_{\mathbb{T}^d} \rho_{in}(x)\theta_\chi^\delta(0,x)dx \quad \forall \chi \in C_c^\infty((0,T) \times \mathbb{T}^d) \quad (4.3)$$

As $\theta_\chi^\delta(0,\cdot)$ converges weakly to $\theta_\chi(0,\cdot)$ in $L^2(\mathbb{T}^d)$, a subsubsequence argument shows that the whole family $\{\tilde{\rho}^\delta\}_{\delta > 0}$ converges weakly-star to $\tilde{\rho}$ in $L^\infty((0,T) \times \mathbb{T}^d)$ as $\delta \downarrow 0$ and

$$\int_0^T \int_{\mathbb{T}^d} \tilde{\rho}(s,x)\chi(s,x)dxds = \int_{\mathbb{T}^d} \rho_{in}(x)\theta_\chi(0,x)dx \quad \forall \chi \in C_c^\infty((0,T) \times \mathbb{T}^d). \quad (4.4)$$

This uniquely characterises $\tilde{\rho}$ independently of the standard mollifier $w$, and in view of (4.2) shows that up to a zero measure set, we have $\rho = \tilde{\rho}$, where we recall that $\rho$ is the vanishing diffusivity solution.

We are missing the proof of no anomalous dissipation. We suppose by contradiction that there exists $\varepsilon > 0$ and a subsequence $v_n \to 0$ such that

$$\limsup_{v_n \to 0} v_n \int_0^T \int_{\mathbb{T}^d} |\nabla \rho^{v_n}(s,x)|^2 dx ds > \varepsilon.$$

In light of Proposition 2.5, the vanishing viscosity solution $\rho$ belongs to $C([0,T]; w - L^2(\mathbb{T}^d))$ and $\rho^v$ converges to $\rho$ in $C([0,T]; w - L^2(\mathbb{T}^d))$ as $v \downarrow 0$. Let us now give a useful lemma.

**Lemma 4.1** *In the context of this section, we have*

$$\lim_{t \to 0^+} \|\rho(t,\cdot)\|_{L^2(\mathbb{T}^d)} = \|\rho_{in}\|_{L^2(\mathbb{T}^d)}.$$

**Proof** We know that $\rho^v \overset{*}{\rightharpoonup} \rho$ in the weak-star topology of $L^\infty((0,T); L^2(\mathbb{T}^d))$ as $v \downarrow 0$. We know $\rho(t,\cdot) \rightharpoonup \rho_{in}$ in the weak $L^2$ topology as $t \downarrow 0$. From the lower semicontinuity of the norms and the bound of Theorem 2.1 we get

$$\|\rho(t,\cdot)\|_{L^2(\mathbb{T}^d)} \leq \|\rho_{in}\|_{L^2(\mathbb{T}^d)} \qquad \forall t \in [0,T].$$

Therefore from the lower semicontinuity of the norms and the previous bound we get

$$\int_{\mathbb{T}^d} |\rho_{in}(x)|^2 dx \geq \limsup_{t \to 0^+} \int_{\mathbb{T}^d} |\rho(t,x)|^2 dx \geq \liminf_{t \to 0^+} \int_{\mathbb{T}^d} |\rho(t,x)|^2 dx \geq \int_{\mathbb{T}^d} |\rho_{in}(x)|^2 dx.$$

□

Firstly, thanks to Lemma 4.1 we fix $\delta > 0$ such that

$$\|\rho(t,\cdot)\|_{L^2(\mathbb{T}^d)}^2 \geq \|\rho_{in}\|_{L^2(\mathbb{T}^d)}^2 - \frac{\varepsilon}{4} \qquad \forall\, t \leq \delta. \tag{4.5}$$

Therefore, we have

$$v_n \int_0^T \int_{\mathbb{T}^d} |\nabla \rho^{v_n}(s,x)|^2 dx ds = v_n \int_0^\delta \int_{\mathbb{T}^d} |\nabla \rho^{v_n}(s,x)|^2 dx ds$$
$$+ v_n \int_\delta^T \int_{\mathbb{T}^d} |\nabla \rho^{v_n}(s,x)|^2 dx ds \tag{4.6}$$

and the first term can be bounded thanks to Theorem 2.4 with

$$2v_n \int_0^\delta \int_{\mathbb{T}^d} |\nabla \rho^{v_n}(s,x)|^2 dx ds = \|\rho_{in}\|_{L^2(\mathbb{T}^d)}^2 - \|\rho^{v_n}(\delta,\cdot)\|_{L^2(\mathbb{T}^d)}^2$$
$$= \|\rho_{in}\|_{L^2(\mathbb{T}^d)}^2 - \|\rho(\delta,\cdot)\|_{L^2(\mathbb{T}^d)}^2$$
$$+ \|\rho(\delta,\cdot)\|_{L^2(\mathbb{T}^d)}^2 - \|\rho^{v_n}(\delta,\cdot)\|_{L^2(\mathbb{T}^d)}^2$$
$$\leq \frac{\varepsilon}{4} + \frac{\varepsilon}{4},$$

where in the last we used (4.5) to estimate $\|\rho(\delta,\cdot)\|_{L^2(\mathbb{T}^d)}^2 - \|\rho^{v_n}(\delta,\cdot)\|_{L^2(\mathbb{T}^d)}^2 \leq \varepsilon/4$ we used the lower semicontinuity of the norms, and have taken $n$ sufficiently large.

To estimate the second term in (4.6) we claim that

$$\limsup_{v_n \to 0} v_n \int_\delta^T \int_{\mathbb{T}^d} |\nabla \rho^{v_n}|^2 dx ds \leq \frac{\varepsilon}{4},$$

which will conclude the proof. We estimate

$$2\nu_n \int_\delta^T \int_{\mathbb{T}^d} |\nabla \rho^{\nu_n}(s,x)|^2 dx ds \leq 2\nu_n \int_\delta^T \int_{\mathbb{T}^d} |\nabla \overline{\rho}^{\nu_n}(s,x)|^2 dx ds$$
$$+ 2\nu_n \int_\delta^T \int_{\mathbb{T}^d} |\nabla(\rho^{\nu_n}(s,x) - \overline{\rho}^{\nu_n}(s,x))|^2 dx ds,$$

where $\overline{\rho}^{\nu_n} \in L^\infty((\delta,T) \times \mathbb{T}^d) \cap L^2((\delta,T); H^1(\mathbb{T}^d))$ is the unique solution to the advection–diffusion ($\nu$-PDE) on $[\delta,1] \times \mathbb{T}^d$ with initial datum $\rho(\delta,\cdot) \in L^2(\mathbb{T}^d)$ by Theorem 2.4. Thanks again to Theorem 2.4 applied to $\overline{\rho}^\nu$ we have

$$\int_{\mathbb{T}^d} |\overline{\rho}^{\nu_n}(T,x)|^2 dx + 2\nu_n \int_\delta^T \int_{\mathbb{T}^d} |\nabla \overline{\rho}^{\nu_n}(s,x)|^2 dx ds \leq \int_{\mathbb{T}^d} |\overline{\rho}^{\nu_n}(\delta,x)|^2 dx. \quad (4.7)$$

Observe that in view of Proposition 2.5, and by uniqueness of bounded weak solutions of (PDE) over $[\delta,T] \times \mathbb{T}^d$, the sequence $(\bar{\rho}^{\nu_n})_{n\in\mathbb{N}}$ converges in $C([\delta,T]; w - L^2(\mathbb{T}^d))$ to $\rho$ as $n \to +\infty$. Hence,

$$\limsup_{\nu_n \to 0} 2\nu_n \int_\delta^T \int_{\mathbb{T}^d} |\nabla \overline{\rho}^{\nu_n}(s,x)|^2 \leq \|\rho(\delta,\cdot)\|^2_{L^2(\mathbb{T}^d)} - \liminf_{\nu_n \to 0} \|\overline{\rho}^{\nu_n}(T,\cdot)\|^2_{L^2(\mathbb{T}^d)} = 0,$$

where in the last equality we use weak lower semicontinuity of the norm, equation (4.7), and conservation of the $L^2$ norm of $\rho$ on $[\delta,T]$ by Theorem 2.1 to get

$$\int_{\mathbb{T}^d} |\rho(\delta,x)|^2 dx = \int_{\mathbb{T}^d} |\rho(T,x)|^2 dx \leq \liminf_{n\to+\infty} \int_{\mathbb{T}^d} |\overline{\rho}^{\nu_n}(T,x)|^2 dx$$
$$\leq \limsup_{n\to+\infty} \int_{\mathbb{T}^d} |\overline{\rho}^{\nu_n}(T,x)|^2 dx \leq \int_{\mathbb{T}^d} |\rho(\delta,x)|^2 dx. \quad (4.8)$$

Finally, from (2.3) applied to $\rho^{\nu_n} - \overline{\rho}^{\nu_n}$ we have that

$$2\nu_n \int_\delta^T \int_{\mathbb{T}^d} |\nabla(\rho^{\nu_n}(s,x) - \overline{\rho}^{\nu_n}(s,x))|^2 dx ds \leq \|\rho^{\nu_n}(\delta,\cdot) - \overline{\rho}^{\nu_n}(\delta,\cdot)\|^2_{L^2(\mathbb{T}^d)}$$
$$= \|\rho^{\nu_n}(\delta,\cdot)\|^2_{L^2(\mathbb{T}^d)} + \|\rho(\delta,\cdot)\|^2_{L^2(\mathbb{T}^d)} - 2\int_{\mathbb{T}^d} \rho^{\nu_n}(\delta,x)\rho(\delta,x) dx, \quad (4.9)$$

where we used that $\overline{\rho}^{\nu_n}(\delta,\cdot) = \rho(\delta,\cdot)$. Using (2.3) and (4.5) we have

$$\limsup_{n\to+\infty} \int_{\mathbb{T}^d} |\rho^{\nu_n}(t,x)|^2 dx \leq \int_{\mathbb{T}^d} |\rho(0,x)|^2 dx \leq \int_{\mathbb{T}^d} |\rho(\delta,x)|^2 dx + \frac{\varepsilon}{4}. \quad (4.10)$$

Taking the lim $\sup_{n\to\infty}$ on both sides of (4.9), using that $\rho^{\nu_n}(\delta,\cdot)$ converges weakly to $\rho(\delta,\cdot)$ by Proposition 2.5 and using (4.10) we conclude the claim and therefore the proof.

**Acknowledgements** GM is supported by the Swiss State Secretariat for Education, Research and Innovation (SERI) under contract number MB22.00034. JP is thankful to Nikolay Tzvetkov for his support. MS was supported by the Swiss State Secretariat for Education, Research and Innovation (SERI) under contract number MB22.00034 during the first part of this work. The authors are thankful to the anonymous referee for useful comments. The authors are thankful to Luigi De Rosa for pointing out an imprecision in a preliminary proof.

**Funding** Open access funding provided by Gran Sasso Science Institute - GSSI within the CRUI-CARE Agreement.